\documentclass[a4paper,11pt]{article}
\usepackage[left=3cm,right=3cm]{geometry}
\usepackage{amssymb}
\usepackage{amsfonts}
\usepackage{amsmath}
\usepackage{bm}
\usepackage{bbm}
\usepackage{latexsym}
\usepackage{epsfig}
\usepackage{algorithm}
\usepackage[noend]{algpseudocode}
\algrenewcommand\textproc{\texttt}
\usepackage{verbatim}
\usepackage{enumerate}
\usepackage{multicol}
\usepackage{mathrsfs}
\usepackage{arydshln}
\usepackage{ulem}
\usepackage{enumitem}
\usepackage{xspace} 
\usepackage{tikz}
\usepackage{pgfplots}
\usepackage{color}

\newcommand{\Infoset}{{\it Info-set}\xspace} 
\newcommand{\RVset}{{\it Velocity-set}\xspace} 
\newcommand{\Residue}{{\it Residue}\xspace} 
\newcommand{\Revector}{{\it Residue Vectors}\xspace} 
\newcommand{\Zmatrix}{{\it Z-matrix}\xspace} 
\newcommand{\Zmatrices}{{\it Z-matrices}\xspace} 

\newcommand{\bM}{{\bf M}}
\newcommand{\bL}{{\bf L}}
\newcommand{\vecV}{{\bf V}}
\newcommand{\vecF}{{\bf F}}

\newcommand{\vecomega}{{\bm \omega}}

\newcommand{\vecf}{{\bf f}}
\newcommand{\vecg}{{\bf g}}

\newcommand{\veczero}{{\bm 0}}

\newcommand{\vecr}{{\bf r}}
\newcommand{\vecres}{{\bm \epsilon}}

\newcommand{\vecv}{{\bf v}}
\newcommand{\vectau}{{\bm \tau}}

\newcommand{\bigO}{{\mathcal O}}
\newcommand\funding[1]{\protect\\ \hspace*{15.37pt}{\bfseries Funding:} #1}

\newtheorem{theorem}{Theorem}
\newtheorem{definition}[theorem]{Definition}

\numberwithin{equation}{section}

\numberwithin{theorem}{section}

\usepackage{amsopn}

\begin{document}


\title{{Hierarchical Orthogonal Matrix Generation and Matrix-Vector
		Multiplications in Rigid Body Simulations}  
	\thanks{Submitted to the editors DATE.
		\funding{This work was funded by the National Science Foundation
			under grant DMS-1217080 (F. Fang and J. Huang) 
			and ACI-1440396 (B. Zhang), and by the Howard Hughes Medical Institute and 
			National Institute of Health under grants GM31749
			and GM103426 (G. Huber and J. A. McCammon). }}}
\author{
	Fuhui Fang \footnotemark[3] \thanks{Corresponding author (fangf@live.unc.edu). }
	\and
	Jingfang Huang \thanks{Department of Mathematics, University of North
		Carolina at Chapel Hill, Chapel Hill, NC 27599-3250.} 
	\and 
	Gary Huber\thanks{Howard Hughes Medical Institute, University of California at 
		San Diego, La Jolla, CA 92093-0365.}
	\and
	J. Andrew McCammon \footnotemark[4]
	\thanks{Department of Chemistry and Biochemistry and Department of
		Pharmacology, University of California at San Diego, La Jolla,
		CA 92093-0365.} 
	\and
	Bo Zhang\thanks{Center for Research in Extreme Scale Technologies,
		Indiana University, Bloomington, IN, 47404.}
}
\maketitle

\begin{abstract}
In this paper, we apply the hierarchical modeling technique and study
some numerical linear algebra problems arising from the Brownian
dynamics simulations of biomolecular systems where molecules are
modeled as ensembles of rigid bodies. Given a rigid body $p$
  consisting of $n$ beads, the $6 \times 3n$ transformation matrix $Z$
  that maps the force on each bead to $p$'s translational and
  rotational forces (a $6\times 1$ vector), and $V$ the row space of
  $Z$, we show how to explicitly construct the $(3n-6) \times 3n$
  matrix $\tilde{Q}$ consisting of $(3n-6)$ orthonormal basis vectors
  of $V^{\perp}$ (orthogonal complement of $V$) using only $\bigO(n \log
  n)$ operations and storage. For applications where only the
matrix-vector multiplications $\tilde{Q}\vecv$ and $\tilde{Q}^T\vecv$
are needed, we introduce asymptotically optimal $\bigO(n)$ hierarchical
algorithms without explicitly forming $\tilde{Q}$.  Preliminary
numerical results are presented  to demonstrate the performance and
accuracy of the numerical algorithms. 
\end{abstract}



\section{Background and problem statement}
\label{sec:problem}
In the Brownian dynamics simulations of biomolecules with hydrodynamic 
interactions, the complex molecular system is modeled as multiple 
(hundreds or thousands) rigid bodies to reduce the numerical ``stiffness'' 
due to the local chemical bond type interactions between atoms that 
cause very high frequency oscillations and subsequently require extremely 
small step size when marching in time. Instead of a thorough listing of 
existing literature on the Brownian dynamics models and hydrodynamic 
interactions, we focus on the ``shell-bead'' model (see, e.g.,
\cite{carrasco1999hydrodynamic,dlugosz2015toward,wang2013assessing})
that describes the hydrodynamic forces exerted on a protein.

In the shell-bead model, the molecular system is represented by $m$
rigid bodies, where rigid body $j$ is modeled by $n_j$ spherical beads
(often of the same radius) placed on the molecular surface with
complex geometry. The total number of beads is 
$n = \sum_{j=1}^m n_j$. Let $\vecf_k^j$ be the external force applied
on bead $k$ (located at $\vecr_k^j$) of rigid body $j$, then rigid body 
$j$'s {\it resultant force} $\vecF_j$ and {\it torque} $\vectau_j$ are given by 
\begin{equation}
\vecF_j = \sum_{k=1}^{n_j} \vecf_k^j, \quad 
\vectau_j= \sum_{k=1}^{n_j} \vecr_k^j \times \vecf_k^j,
\end{equation}
or in matrix form 
\begin{equation}
\left[ \begin{array}{c} 
  \vecF_j \\ 
  \vectau_j
  \end{array}
\right]_{6 \times 1} 
=  
\left[ \begin{array} {cccc}
I&I& \dots &I \\
A_1^j & A_2^j & \dots &A_{n_j}^j \end{array} \right]_{6 \times 3n_j}  
  \left[ \begin{array}{c}
  \vecf_1^j \\ \vecf_2^j \\ \cdots \\  \vecf_{n_j}^j \end{array} 
    \right]_{3n_j \times 1} 
= Z^j \vecf^j,
\end{equation}
where we refer to the $6 \times 3n_j$ matrix as the \Zmatrix. In the formula, $I$ 
is the $3\times 3$ identity matrix and $A_k^j$ is the transformation matrix
of $A_k^j \cdot \vecf_k^j \triangleq \vecr_k^j \times\vecf_k^j$. 
Particularly, if $\vecr_k^j = (x_k^j, y_k^j, z_k^j)$, then 
\[
A_k^j = \begin{bmatrix} 
0 & -z_k^j & y_k^j \\
z_k^j & 0 & -x_k^j \\
-y_k^j & x_k^j & 0 
\end{bmatrix}.
\]
The shell-bead model assumes that the hydrodynamic effects are related
to the deterministic forces through $D\vecf = \vecv$, where the 
$3n \times 3n$ matrix $D$ is the symmetric Rotne-Prager-Yamakawa
tensor whose entries are determined by the bead locations (see, e.g., 
\cite{batchelor1976brownian,ermak1978brownian}), and $\vecf$ and
$\vecv$ in Matlab notation are given by 
\[
\vecf  = [\vecf_1^1; \cdots; \vecf_{n_1}^1; \cdots
         \vecf_1^m; \cdots; \vecf_{n_m}^m], \quad 
\vecv  = [\vecv_1^1; \cdots; \vecv_{n_1}^1; \cdots; 
         \vecv_1^m; \cdots; \vecv_{n_m}^m]. 
\]
Let $\Big [ [\vecF_1; \vectau_1]; \cdots; [\vecF_m; \vectau_m] \Big]$
be the external deterministic force and torque vector acting on the
$m$ rigid bodies. Under the rigid body constraint, the corresponding
deterministic velocity field of the rigid body can be obtained by
solving 
\begin{equation}
\label{eq:rigidbodyBDHI}
\underbrace{
  Z_{6m \times 3n} 
  \underbrace{ 
    (D^{-1})_{3n \times 3n} 
    \underbrace{
     (Z^T)_{3n \times 6m} 
     \begin{bmatrix} 
     \vecV_1 \\
     \vecomega_1 \\
     \vdots \\
     \vecV_m \\
     \vecomega_m 
     \end{bmatrix} 
    }_{\mbox{velocities on beads}}    
  }_{\mbox{forces on beads}}
}_{\mbox{forces on rigid bodies}} = 
  \begin{bmatrix} 
    \vecF_1 \\
    \vectau_1 \\
    \vdots \\
    \vecF_m \\
    \vectau_m 
  \end{bmatrix}, 
\end{equation}
where $Z = {\rm diag}(Z^1, \cdots, Z^m)$ is block diagonal. Here, the force and
translational velocity field are denoted by uppercase letters for the
rigid bodies and lowercase letters for the beads. Eq.~(\ref{eq:rigidbodyBDHI})
simply means that the unknown rigid body velocity field, when mapped
onto individual beads, should yield the force field on the beads via
$\vecf = D^{-1}\vecv$. Then, the force acting on the rigid body can be
obtained by integrating all the bead forces using $Z\vecf$ and the
result should match the given external forces on the right-hand side. 

One major numerical difficulty in solving Eq.~(\ref{eq:rigidbodyBDHI})
accurately and efficiently is the calculation of $D^{-1}$ as $D$ is
dense. For large $n$,  even with the acceleration of the 
fast direct solvers~\cite{greengard, ho2012fast} or 
$H$-matrix techniques~\cite{hackbusch1999sparse,hackbusch2000sparse}, 
computing $D^{-1}$ at each time marching step is simply too expensive 
for dynamic simulations. It is possible to avoid the explicit computation 
of $D^{-1}$ by a reformulation of Eq.~(\ref{eq:rigidbodyBDHI}). To simplify the 
discussion, we assume that with an easy QR procedure on both sides of 
Eq.~(\ref{eq:rigidbodyBDHI}) for each rigid body $j$, the $6$ vectors in the 
diagonal block $Z^j$ become orthogonal. Denote the orthogonal version of the
matrix $Z^j$ by $Q_Z^j$. Define $Q_Z={\rm diag}(Q_Z^1, \cdots, Q_Z^m)$ and 
$\tilde{Q}$ the orthogonal vectors such that $Q = [Q_Z; \tilde{Q}]$ is a 
$3n \times 3n$ orthogonal matrix. Here, we mildly abuse the notations by using 
$\vecV = [\vecV_1; \vecomega_1; \cdots; \vecV_m; \vecomega_m]$ and introducing $\vecF$ 
to represent the new right hand side after the QR process. Pad $(3n-6m)$ zeros to $\vecV$ 
and introduce $(3n-6m)$ unknowns $\vecg$ to $\vecF$, we form 
\begin{equation}
  \begin{bmatrix} 
    Q_Z \\ \tilde{Q} 
  \end{bmatrix} 
  D^{-1} 
  \begin{bmatrix} 
   Q_Z^T & \tilde{Q}^T 
  \end{bmatrix} 
  \begin{bmatrix} 
    \vecV\\
    \veczero 
  \end{bmatrix} 
 = \begin{bmatrix} 
    \vecF \\
    \vecg 
   \end{bmatrix}.
\end{equation}
As the inverse of the orthogonal matrix $Q$ is simply its transpose,
some algebraic manipulations show that for any given vector $\vecF$,
one can first find the unknown $\vecg$ by solving 
\begin{equation}
  \label{eq:newBDHI}
  \veczero = \tilde{Q} D Q_Z^T \vecF + \tilde{Q}D\tilde{Q}^T \vecg
\end{equation}
via a preconditioned Krylov subspace iterative method. Physically, the 
vector $\vecg$ can be envisioned as the constraint forces needed to
keep the beads in each rigid body together during the imposed external forces.
Then, the velocity vector $\vecV$ can be computed by 
\begin{equation}
  \label{eq:solution}
   \vecV = Q_Z D(Q_Z^T \vecF + \tilde{Q}^T \vecg).
\end{equation}

The fundamental building blocks required by this new formulation are the 
fast matrix-vector multiplications of $D\vecv$, $\tilde{Q}^T\vecv$, and 
$\tilde{Q}\vecv$ with any given vector $\vecv$. The $D\vecv$ operation can 
be carried out efficiently using the fast multipole methods as discussed in 
\cite{liang2013fast,Ying2004}. This paper considers the $\tilde{Q} \vecv$ 
and $\tilde{Q}^T \vecv$ operations. Notice that the matrix $Z$ is block 
diagonal and each block corresponds to a rigid body, this structure
allows us to construct $\tilde{Q}^j$ for each $Q_Z^j$ separately and
then form the $\tilde{Q}$ matrix using 
$\tilde{Q}={\rm diag}(\tilde{Q}^1, \cdots, \tilde{Q}^m)$. We therefore 
focus on one rigid body and drop the index $j$ in the following problem
statement.

{\bf \noindent Problem Statement:} Given the \Zmatrix of a rigid body $p$ 
with $n$ beads, how to efficiently construct and store $\tilde{Q}$ 
explicitly (if needed)? And for any given vector $\vecv$ of proper size, 
how to efficiently compute $\tilde{Q} \vecv$ and $\tilde{Q}^T \vecv$? 

The main contributions of this paper are three novel algorithms optimal in
complexity and storage requirements. The discussions of these algorithms are 
organized as follows. In Sec.~2, we introduce the hierarchical tree structure 
and present a hierarchical model for constructing $\tilde{Q}$. To better preserve
orthogonality, in Sec.~3, we apply tools from the orthogonal linear algebra and present 
the first algorithm to explicitly construct $\tilde{Q}$ using $\bigO(n \log n)$ operations 
and storage. In Sec.~4, we present two asymptotically optimal algorithms to compute the 
matrix-vector multiplications $\tilde{Q}\vecv$ and $\tilde{Q}^T\vecv$ using only 
$\bigO(n)$ operations and storage, without explicitly forming $\tilde{Q}$. In Sec.~5,
numerical results are presented to demonstrate the algorithms' performance and 
orthogonality properties. Finally in Sec.~6, we summarize our results and discuss several related
research topics.

\section{Hierarchical tree and hierarchical model}
We apply the hierarchical modeling technique to study the orthonormal
basis vectors in $\tilde{Q}$.  The hierarchical modeling technique
identifies any low-rank, or low-dimensional, or other compact features
in the system, and the compressed representations are then recursively
collected from children to parents, and transmitted between different
nodes on a hierarchical tree structure using properly compressed
translation operators. In its numerical implementation, the
hierarchical models are often re-expressed as recursive algorithms,
which can be easily interfaced with existing dynamical schedulers from
High-Performance Computing (HPC) community for optimal parallel
efficiency.

Different aspects of the hierarchical modeling technique have been
known and addressed by different research communities
previously. Examples include the classical fast Fourier transform
(FFT) \cite{cooley1965algorithm} where the Halving Lemma shows how
data can be compressed and the odd-even term splitting of the
polynomials creates a hierarchical tree to allow recursively
processing the compressed information efficiently; the multigrid
method (MG) \cite{brandt1977multi,hackbusch2013multi} where the
hierarchical tree structure is formed via adaptively refining the
computational domain, and data compression and transmission are
performed using the relaxation (smoother) and projection (restriction)
operators by analyzing the frequency domain behaviors of the error
functions between different levels of the (adaptive) tree to
effectively reduce the high frequency errors; and the fast multipole
method (FMM) \cite{greengard1987fast,greengard1997new} where the
particle information inside a box is first compressed to the multipole
expansion, and the compressed information is transmitted recursively
to parent levels using the multipole-to-multipole translation in the
upward pass on the hierarchical tree structures. The collected and
compressed information is later transmitted to target boxes using the
multipole-to-local translations and propagated to child levels using
the local-to-local translations in the downward pass.  When there are
$n$ terms (FFT) in the polynomial or $n$ approximately uniformly
distributed particles (MG or FMM), the depth of the hierarchical tree
is normally $\bigO(\log n)$ and the number of tree nodes is
approximately $\bigO(n)$.  Therefore, if each level only requires
$\bigO(n)$ operations (e.g., FFT), the algorithm complexity will be
$\bigO(n \log n)$. If each tree node only requires a constant amount
of operations (e.g. MG or FMM), the algorithm complexity will be
asymptotically optimal $\bigO(n)$. In this section, we discuss how to
use the hierarchical modeling technique to answer the questions in
the Problem Statement in Sec.~\ref{sec:problem}.

\subsection{Adaptive hierarchical tree structure}
We first consider generating a spatial adaptive hierarchical tree 
when simulating a molecular system modeled by multiple rigid bodies in
the shell-bead model. We assume each rigid body is ``discretized" into
a number of beads to capture the hydrodynamic interactions between
rigid bodies.  A hierarchical partition is then performed to divide
the beads domain into nested cubical boxes, where the root box is the
smallest bounding box that contains all the beads. Without loss of
generality, the root box is normalized to size $1$ in each side. The
root box is partitioned equally along each dimension. The partition
continues recursively on the resulting box until the box 
contains no more than $s$ beads,
at which point it becomes a leaf node. Empty boxes encountered during
partition are pruned off. Here, we set $s=1$ to simplify the
discussions in the following sections. In our implementation, other
values of $s$ are allowed after modifying how the leaf nodes are
processed. 

{\noindent \bf Comment:} The octree can be modified to form a binary
tree. At a parent node $p$, one can first create two ``ghost'' nodes
by separating the beads within $p$ by the $z$-direction. Then, each
ghost node can be separated by the $y$-direction,
creating four more ghost nodes. Finally, these four ghost nodes are
partitioned along the $x$-direction, creating the actual eight child
nodes of $p$. The depth of the binary tree is at most $3L$, where 
$L = O(\log n)$. Such modification is not fundamental, but could significantly
simplify both the notations and descriptions of the algorithms in the
rest of the paper. For this reason, we focus on the binary tree and
set $s=1$.

\subsection{Divide-and-conquer strategy and hierarchical model}
\label{sec:divideconquer}
We consider a particular choice of the orthonormal vectors in
$\tilde{Q}$ using the divide-and-conquer strategy on the {\it
  hierarchical} tree structure. We start from a two level setting
where the parent rigid body $p$ consisting of $n$ beads is partitioned
into two child nodes, child $x$ with $n_x$ beads and child $y$ with
$n_y$ beads, where $n_x + n_y = n$. Let 
\[
Z_x = \begin{bmatrix} 
I & I & \cdots & I \\
A_1 & A_2 & \cdots & A_{n_x} 
\end{bmatrix}, \quad 
Z_y = \begin{bmatrix} 
I & I & \cdots & I \\
B_1 & B_2 & \cdots & B_{n_y} 
\end{bmatrix}
\]
be the \Zmatrices of $x$ and $y$, respectively. Assume both $Z_x$ and
$Z_y$ are full rank, the orthogonal matrices $\tilde{Q}_x$ of size
$(3n_x - 6) \times 3n_x$ and $\tilde{Q}_y$ of size 
$(3n_y - 6) \times 3n_y$ that satisfy $\tilde{Q}_x \perp Z_x$ and 
$\tilde{Q}_y \perp Z_y$ are available in compact form. The key
observation comes from the study of the matrix 
\[
H = \begin{bmatrix} 
Z_x & Z_y \\
\veczero & Z_y \\
\tilde{Q}_x & \veczero \\
\veczero & \tilde{Q}_y 
\end{bmatrix}_{3n \times 3n}
\]
and the fact that $Z_p=[Z_x,Z_y]$. 
It is straightforward to verify that the vectors in the lower
$(3n-12)$ rows of $H$ are normalized, orthogonal to each other and to the
first 12 rows of $H$. This means that $\tilde{Q}_p$ of  parent $p$ can
readily ``receive'' the lower $(3n-12)$ rows of vectors from its two
children. For the remaining $6$ row vectors in $\tilde{Q}_p$, a
Gram-Schmidt procedure on the first 12 row vectors can be performed
and the last $6$ orthonormal vectors will be orthogonal to the vectors
in $Z_p$ and the $(3n-12)$ vectors from the children. In a multilevel
setting for a rigid body with $n$ beads, there will be approximately
$\bigO(\log n)$ levels and for each level, the Gram-Schmidt procedure
requires approximately $\bigO(n)$ operations and storage to explicitly
generate all the orthogonal vectors for that level. The total storage
and operations required are therefore both $\bigO(n \log n)$. 

Unfortunately,  straightforward implementation of this
 divide-and-conquer idea will result in an
algorithm with stability issues. In particular, the matrix $\tilde{Q}$
will lose orthogonality and $\tilde{Q} \tilde{Q}^T \neq I$. 
One source for the instability is the ill-conditioning
of the \Zmatrices. For instance, when all the beads are located 
exactly on a straight line, 
the rank of the \Zmatrix \ is only $5$ instead of $6$, 
and the $H$ matrix becomes singular. In the 
following sections, we show how to resolve the instability issues
using the orthogonal linear algebra techniques. 

\section{Stable $\bigO(n \log n)$ orthogonal matrix generation
  algorithm}
We start from Theorem~\ref{thm:masscenter}. 
\begin{theorem}
\label{thm:masscenter}
Assume the centroid of a rigid body $p$ is located at the origin
$\veczero$, then the first $3$ rows of the \Zmatrix are orthogonal
to the last $3$ rows. 
\end{theorem}
The proof follows from the identifies 
\[
\sum x_i = \sum y_i = \sum z_i = 0 
\]
if the centroid of the rigid body is chosen as the origin. It suggests
that in order to get $Q_Z$, the orthogonal version of the \Zmatrix of
$p$,  the Gram-Schmidt procedure only needs to be applied to the last
$3$ rows of $Z_p$. We adopt this assumption in the following
discussions and formulas and revisit the two level setting in 
Sec.~\ref{sec:divideconquer}. To preserve orthogonality, we assume 
$Z_x$ and $Z_y$ are already decomposed with respect to their
respective centers as 
\begin{equation}
\label{eq:Zxdecomp}
	Z_x = C_x\cdot Q_{Z_x}
= \left[
  \begin{array}{c:c} 
    \sqrt{n_x} I_{3 \times 3} & \veczero \\ 
    \hdashline \veczero & R_{22} 
  \end{array} 
  \right] 
  \cdot \left[ 
    \begin{array}{c} I_{n_x} \\ 
      \hdashline V_{3 \times 3n_x} \end{array} 
    \right]
\end{equation}
and 
\begin{equation}
	\label{eq:Zydecomp}
	Z_y = C_y \cdot Q_{Z_y} 
  = \left[
    \begin{array}{c:c} 
      \sqrt{n_y} I_{3 \times 3} & \veczero \\ 
      \hdashline \veczero & S_{22} 
  \end{array} 
    \right]
  \cdot \left[ 
    \begin{array}{c} I_{n_y}  \\ 
      \hdashline W_{3 \times 3n_y} \end{array} 
    \right],
\end{equation}
where
\[
I_{n_x} = \frac{1}{\sqrt{n_x}} 
\begin{bmatrix} I_{3\times 3} & \cdots & I_{3 \times 3} \end{bmatrix}_{3\times 3n_x}, 
\quad 
I_{n_y} = \frac{1}{\sqrt{n_y}} 
\begin{bmatrix} I_{3\times 3} & \cdots & I_{3 \times 3} \end{bmatrix}_{3\times 3n_y}.
\]
In the formula, as each child uses its own centroid, the first $3$
rows of the \Zmatrix are orthogonal to the last $3$ rows by
Theorem~\ref{thm:masscenter} and matrix $C$ is block diagonal. 

Assuming orthogonality preserving results are available for $x$ and
$y$, i.e., the orthogonal submatrices 
\[
Q_{Z_x} = \begin{bmatrix} I_{n_x} \\ V \end{bmatrix}_{6 \times 3 n_x},
\quad 
\tilde{Q}_x \perp Q_{Z_x}, \quad 
Q_{Z_y} = \begin{bmatrix} I_{n_y} \\ W \end{bmatrix}_{6 \times 3n_y}, 
\quad 
\tilde{Q}_y \perp Q_{Z_y}
\]
are already constructed with excellent orthogonality properties, we
study in the following how to stably find $\tilde{Q}_p$ as well as
$T_{22}$ and $Q_{Z_p}$ in the decomposition of $p$'s \Zmatrix
\begin{equation} 
\label{eq:zporthdecomp}
  Z_p = C_p\cdot Q_{Z_p} = 
  \left[
    \begin{array}{c:c} 
      \sqrt{n} I_{3 \times 3} & \veczero \\ 
      \hdashline
        \veczero &T_{22} \end{array} \right]
  \cdot \left[ 
	\begin{array}{c} I_{n_p} \\ 
      \hdashline U_{3 \times 3n}
        \end{array} 
    \right], 
	\end{equation}
	\begin{equation*}
	\quad I_{n_p} = \frac{1}{\sqrt{n}} 
\begin{bmatrix} 
I_{3 \times 3} & \cdots & I_{3 \times3} 
\end{bmatrix}_{3\times 3n},
\end{equation*} 

where the new origin is located at the centroid of $p$ that can be
easily computed from the centroids of $x$ and $y$. Following the ideas
in  Sec.~\ref{sec:divideconquer}, we consider
a particular orthogonal matrix $Q = [Q_{Z_p}; \tilde{Q}_p]$, where 
\[
\tilde{Q}_p = \begin{bmatrix} 
  \Revector \\
  \begin{array}{cc} 
    \tilde{Q}_x & \veczero \\
    \veczero & \tilde{Q}_y 
  \end{array} 
\end{bmatrix} 
\]
and the row vectors in $[Q_{Z_p}; \Revector]$ form the same subspace
as that spanned by the 12 row vectors in 
$
\left[ \begin{array}{cc} Q_{Z_x}  & \veczero \\ 
    \veczero & Q_{Z_y} \end{array} \right].
$ 
Clearly, $\tilde{Q}_p$ contains two parts, the lower $3n-12$ rows 
are processed at child levels and require no additional operations or
storage (only nonzero values are stored), and the first $6$ row
vectors are the ``left-over" vectors (referred to as the \Revector) 
after identifying and removing $Q_{Z_p}$ components from the subspace of 
dimension $12$. Therefore, in the hierarchical modeling technique, 
given compressed information $Q_{Z_x}$, $Q_{Z_y}$, $R_{22}$ and
$S_{22}$ at the child level, we study how to compute parent's compressed 
information $Q_{Z_p}$, $T_{22}$ as well as the one time output \Revector.

We prefer the well-conditioned $Q_{Z_x}$ and $ Q_{Z_y}$ to preserve
orthogonality properties in the hierarchical model to the 
original $Z_x$ and $Z_y$ discussed in Sec.~\ref{sec:divideconquer}
because the latter may become rank deficient if all the beads are located
on the same line. Note that the subspace spanned by the \Zmatrix  
is always a subset of the subspace spanned by the corresponding
$Q_Z$. Also, to take advantage of Theorem~\ref{thm:masscenter}, 
instead of considering the row vectors in 
$
\left[ \begin{array}{cc} Q_{Z_x}  & \veczero  \\ 
\veczero & Q_{Z_y} \end{array} \right]
$, 
we consider the orthogonal vectors in the matrix $C$ presented in 
Theorem~\ref{thm:orthovector}, where the 
first $3$ row vectors are simply $I_{n_p}$. 
This strategy saves operations when performing the 
QR decomposition.
\begin{theorem}
\label{thm:orthovector}
For a parent node $p$ with child nodes $x$ and $y$, the 12 row vectors 
of the matrix
\[
C = \left[ 
  \begin{array}{cc}
  \sqrt{{n_x}/{n}} I_{n_x}  & \sqrt{{n_y}/{n}} I_{n_y}  \\
  V & \veczero \\
  \veczero & W \\
  \sqrt{{n_y}/{n}} I_{n_x}  & -\sqrt{{n_x}/{n}} I_{n_y}  \\
\end{array}
\right]
\]
are orthonormal, span the same subspace as the row vectors in 
$
\left[ 
\begin{array}{cc} 
  Q_{Z_x}  & \veczero \\ 
  \veczero & Q_{Z_y} 
\end{array} 
\right]
$,
and the first $3$ rows satisfy 
	$\left[\sqrt{{n_x}/{n}} I_{n_x} ,  \sqrt{{n_y}/{n}}  I_{n_y} \right] = I_{n_p}$.
\end{theorem}
This theorem is simply the result of $C\cdot C^T=I_{12\times 12}$. 
As the orthonormal basis vectors in $C$ are constructed analytically, 
their orthogonality properties are well-preserved.

To compute $Q_{Z_p}$, $T_{22}$, and the \Revector, we first represent
$Z_p$ using $Z_x$ and $Z_y$ by shifting the centroid using the $3
\times 3$ matrices $R_{21}$ and $S_{21}$ as in 
\[
Z_p = \begin{bmatrix} 
\begin{bmatrix} I & \veczero \\ R_{21} & I \end{bmatrix} \cdot Z_x, 
\begin{bmatrix} I & \veczero \\ S_{21} & I \end{bmatrix} \cdot Z_y
\end{bmatrix}.  
\]
Substitute the orthogonal decompositions of $Z_x$ and $Z_y$ available
in Eqs.~(\ref{eq:Zxdecomp}) and (\ref{eq:Zydecomp}), we have
\begin{equation}
\label{eq:Zp}
Z_p 
= \left[ 
  \begin{array}{c:c}
    \sqrt{n_x}I_{n_x} & \sqrt{n_y}I_{n_y}\\ 
    \hdashline
    \sqrt{n_x} R_{21} I_{n_x}+R_{22}V & 
    \sqrt{n_y} S_{21}I_{n_y}+S_{22}W 
  \end{array}
  \right].
\end{equation}
Notice that $I_{n_x} \perp V$, $I_{n_y} \perp W$, and 
\[
\begin{bmatrix} \sqrt{n_x}I_{n_x} & \sqrt{n_y}I_{n_y} \end{bmatrix} 
\perp
\begin{bmatrix} 
\sqrt{n_x} R_{21}I_{n_x} + R_{22}V & 
\sqrt{n_y} S_{21}I_{n_y} + S_{22}W 
\end{bmatrix}, 
\]
One can easily derive 
\[
n_x R_{21}I_{n_x} + n_y S_{21} I_{n_y} = \veczero \Rightarrow 
R_{21} = -\frac{n_y}{n_x} S_{21}. 
\]
As the first $3$ rows of $Z_p(1:3,:)$ are simply $\sqrt{n} I_{n_p}$, 
by Theorem~\ref{thm:masscenter}, we only need to consider the last
$3$ rows of $Z_p$ in Eq.~(\ref{eq:Zp}) reformulated as 
\[
Z_p(4:6,:) = \begin{bmatrix} 
R_{22} & S_{22} & \sqrt{\frac{n_x n}{n_y}} R_{21} 
\end{bmatrix} 
\begin{bmatrix} 
V & \veczero \\ \veczero & W \\ 
\sqrt{\frac{n_y}{n}} I_{n_x} & 
-\sqrt{\frac{n_x}{n}} I_{n_y} 
\end{bmatrix}
\]
using the orthonormal vectors in the $C$ matrix in Theorem~\ref{thm:orthovector}.
Applying existing orthogonality preserving QR algorithms from the
orthogonal linear algebra packages (e.g., Matlab {\tt qr} command), 
we can derive the QR decomposition of the $3\times 9$ matrix 
\begin{equation}
\label{eq:Tformula}
\begin{bmatrix} 
R_{22} & S_{22} & \sqrt{\frac{n_x n}{n_y} } R_{21}
\end{bmatrix} = 
\begin{bmatrix} 
  \begin{array}{c:c} 
    \begin{array}{ccc}
      r_{11}&r_{21}&r_{31}\\
      r_{12}&r_{22}&r_{32}\\
      r_{13}&r_{23}&r_{33}
    \end{array} & 
    \veczero  
  \end{array}
 \end{bmatrix}_{3 \times 9} 
\begin{bmatrix} 
    Q_{11}^T & Q_{12}^T & Q_{13}^T\\
    Q_{21}^T & Q_{22}^T & Q_{23}^T\\
    Q_{31}^T & Q_{32}^T & Q_{33}^T
\end{bmatrix}_{9 \times 9} 
\end{equation}
where each $Q_{ij}$ is a $3 \times 3$ matrix. Then, the last $3$
rows of $Z_p$ become 
\begin{align*} 
	& Z_p(4:6,:)= \begin{bmatrix} 
    \sqrt{n_x}R_{21} I_{n_x} + R_{22} V & 
    \sqrt{n_y}S_{21} I_{n_y} + S_{22} W 
  \end{bmatrix} \\
= & \begin{bmatrix} 
  \begin{array}{c:c} 
    \begin{array}{ccc}
      r_{11}&r_{21}&r_{31}\\
      r_{12}&r_{22}&r_{32}\\
      r_{13}&r_{23}&r_{33}
    \end{array} & 
    \veczero  
  \end{array}
 \end{bmatrix}
\begin{bmatrix} 
Q_{11}^T V + \sqrt{\frac{n_y}{n}}Q_{13}^TI_{n_x} & 
Q_{12}^T W - \sqrt{\frac{n_x}{n}}Q_{13}^TI_{n_y} \\
Q_{21}^T V + \sqrt{\frac{n_y}{n}} Q_{23}^T I_{n_x} & 
Q_{22}^T W - \sqrt{\frac{n_x}{n}} Q_{23}^T I_{n_y} \\
Q_{31}^T V + \sqrt{\frac{n_y}{n}} Q_{33}^T I_{n_x} &
Q_{32}^T W - \sqrt{\frac{n_x}{n}} Q_{33}^T I_{n_y} 
\end{bmatrix}.
\end{align*} 
Notice that the row vectors in the second matrix are orthonormal and
comparing with Eq.~(\ref{eq:zporthdecomp}), we derive 
\begin{align} 
T_{22} & = \begin{bmatrix} 
r_{11} & r_{21} & r_{31} \\
r_{12} & r_{22} & r_{32} \\
r_{13} & r_{23} & r_{33} 
\end{bmatrix},
\label{eq:T22full} \\
U & = \begin{bmatrix} 
Q_{11}^TV + \sqrt{\frac{n_y}{n}}Q_{13}^TI_{n_x} & 
Q_{12}^TW - \sqrt{\frac{n_x}{n}}Q_{13}^TI_{n_y} 
\end{bmatrix}, 
\label{eq:Ufull}\\
\Revector & = \begin{bmatrix} 
Q_{21}^T V + \sqrt{\frac{n_y}{n}} Q_{23}^T I_{n_x} &
 Q_{22}^T W -\sqrt{\frac{n_x}{n}} Q_{23}^T I_{n_y} \\
Q_{31}^T V + \sqrt{\frac{n_y}{n}} Q_{33}^T I_{n_x} &
 Q_{32}^T W -\sqrt{\frac{n_x}{n}} Q_{33}^T I_{n_y}
\end{bmatrix} .
\label{eq:residuefull}
\end{align} 

There are two cases that require special treatments in the
algorithm.  The first case is when both $x$ and $y$ are leaf nodes,
containing one bead each. The second case is $x$ is a leaf node with a
single bead and $y$ is a multi-bead rigid body, or vice versa. The
special treatments are presented as follows. \\

{\noindent \bf Case I: Bead \& Bead.} Consider a parent node with two
child leaf nodes $x$ and $y$, containing a single bead each. The
locations of the beads are given by $(a, b, c)$ and $(-a, -b, -c)$
such that the parent's centroid is located at the origin. Assume
further that $c \neq 0$, then the parent's $U$ and $T_{22}$ matrices
are given by 
\begin{align}
T_{22} &= \begin{bmatrix} 
\sqrt{2(b^2+c^2)} & 0 & 0 \\
-\frac{\sqrt{2}ab}{\sqrt{b^2+c^2}} & 
\sqrt{\frac{2c^2(a^2+b^2+c^2)}{b^2+c^2}} & 0\\
-\frac{\sqrt{2}ac}{\sqrt{b^2+c^2}} &
-\sqrt{\frac{2b^2(a^2+b^2+c^2)}{b^2+c^2}} & 0 \end{bmatrix},
\label{eq:T22beadbead}\\
U & = \begin{bmatrix} 
  \frac{1}{\sqrt{2(b^2+c^2)}} 
  \begin{bmatrix} 0 & -c & b & 0 & c & -b \end{bmatrix} \\
  \frac{1}{\sqrt{2(a^2+b^2+c^2)(b^2+c^2)}} 
  \begin{bmatrix}  b^2+c^2 & -ab & -ac & -b^2-c^2 & ab & ac \end{bmatrix} \\
  \frac{1}{\sqrt{2(a^2+b^2+c^2)}}
  \begin{bmatrix} -a & -b & -c & a & b & c \end{bmatrix}
  \end{bmatrix}, 
\label{eq:Ubeadbead}
\end{align}
and there is no \Revector  generated. Note that in this
case the rank of $T_{22}$  is only two, but $U$ is orthogonal with full rank.\\

{\noindent \bf Case II: Rigid Body \& Bead.} 
Assume child $x$ has $n_x > 1$ beads and child $y$ has a single bead,
$Z_p$ can be reformulated as 
\[
\left[\begin{array}{c:c}
 \sqrt{n_x} I_{n_x} & I \\ \hdashline
 \sqrt{n_x} R_{21} I_{n_x}+R_{22}V & S_{21} 
  \end{array} 
  \right].
\]
Applying QR algorithm from the orthogonal linear algebra package, the
last $3$ rows of $Z_p$ could be reformulated as 
\begin{align*} 
&Z_p(4:6,:) =  \left[ 
    \begin{array}{c:c} \sqrt{n_x} R_{21} I_{n_x} + R_{22}V &
      S_{21} \end{array}
  \right] \\ 
= & 
\begin{bmatrix} R_{22} & \sqrt{n_x n} R_{21} \end{bmatrix} 
\begin{bmatrix} 
V & \veczero \\ 
\sqrt{\frac{1}{n}} I_{n_x} & -\sqrt{\frac{n_x}{n}} I 
\end{bmatrix} \\
= & 
\left[\begin{array}{c:c}
\begin{array}{ccc}
r_{11}&r_{21}&r_{31}\\
r_{12}&r_{22}&r_{32}\\
r_{13}&r_{23}&r_{33}
\end{array} & \veczero  
\end{array}\right]_{3\times 6} 
\begin{bmatrix} 
Q_{11}^T & Q_{12}^T \\ 
Q_{21}^T & Q_{22}^T 
\end{bmatrix}_{6 \times 6} 
\begin{bmatrix} 
  V & \veczero  \\ 
\sqrt{\frac{1}{n}} I_{n_x} & - \sqrt{\frac{n_x}{n}} I 
\end{bmatrix}_{6 \times 6n} \\
= & \left[\begin{array}{c:c}
\begin{array}{ccc}
r_{11}&r_{21}&r_{31}\\
r_{12}&r_{22}&r_{32}\\
r_{13}&r_{23}&r_{33}\end{array}&\veczero  
\end{array}\right]
\begin{bmatrix} 
  Q_{11}^T V + \sqrt{ \frac{1}{n}} Q_{12}^T I_{n_x} & 
 - \sqrt{ \frac{n_x}{n}} Q_{12}^T \\
  Q_{21}^T V + \sqrt{ \frac{1}{n}} Q_{22}^T I_{n_x} & 
 - \sqrt{ \frac{n_x}{n}} Q_{22}^T 
\end{bmatrix}. 
\end{align*}
The orthogonal decomposition and \Revector of parent $p$ are given by 
\begin{align}
T_{22} = &  \begin{bmatrix} 
r_{11}&r_{21}&r_{31}\\
r_{12}&r_{22}&r_{32}\\
r_{13}&r_{23}&r_{33}
\end{bmatrix}, 
\label{eq:T22beadrigid}\\
U = & \begin{bmatrix} 
 Q_{11}^T V + \sqrt{ \frac{1}{n}} Q_{12}^T I_{n_x} & 
 - \sqrt{ \frac{n_x}{n}} Q_{12}^T \end{bmatrix}_{3 \times 3n}
	\label{eq:Ubeadrigid}, \mbox{ and } \\
\Revector = & \begin{bmatrix} 
Q_{21}^T V + \sqrt{ \frac{1}{n}} Q_{22}^T I_{n_x} & 
- \sqrt{ \frac{n_x}{n}} Q_{22}^T \end{bmatrix}_{3 \times 3n}.
\label{eq:residuebeadrigid}
\end{align} 
Note that the dimension of the \Revector in this case is $3$.

Given an adaptive binary tree structure, Algorithm~\ref{alg:explicitQ}
shows the pseudocode of the recursive function {\tt Q\_gen} that
explicitly generates the matrix $Q = [Q_Z; \tilde{Q}]$. The orthogonal
matrix is generated by calling {\tt Q\_gen} on the root node. 

\renewcommand{\algorithmiccomment}[1]{\hfill//\textit{#1}}
\begin{algorithm}[ht]
\begin{algorithmic}[1]
\Function{Q\_gen}{$p$}
\If{$p$ is leaf node} 
\State centroid = bead location, $U = \veczero$,  $T_{22} = \veczero_{3\times3}$
\Else 
\State Find child nodes $x$ and $y$ of node $p$
\State \Call{Q\_gen}{$x$}
\State \Call{Q\_gen}{$y$}
\State Compute centroid of $p$ and form $R_{21}$ and $S_{21}$
\If{both $x$ and $y$ are leaf nodes}
\State Compute $T_{22}$ and $U$ using Eqs.~(\ref{eq:T22beadbead}) and
(\ref{eq:Ubeadbead}) 
\ElsIf {only one of $x$ and $y$ is a leaf node}
\State Compute $T_{22}$, $U$, and \Revector using
Eqs.~(\ref{eq:T22beadrigid}), (\ref{eq:Ubeadrigid}),  
  and (\ref{eq:residuebeadrigid})
\State Output \Revector of size $3$
\Else 
\State Compute $T_{22}$, $U$, and \Revector   using
Eqs.~(\ref{eq:T22full}), (\ref{eq:Ufull}) and (\ref{eq:residuefull}) 
\State Output \Revector of size $6$ 
\EndIf
\EndIf
\If{$p$ is the root} 
\State Output $U$

\EndIf
\EndFunction
\end{algorithmic}
\caption{Recursive algorithm for explicit generation of $Q$}
\label{alg:explicitQ} 
\end{algorithm}

{\noindent \bf Algorithm Complexity.} 
 To estimate the algorithm complexity and storage requirement, 
we consider a system with $n$ beads and a tree with $\bigO(\log n)$
levels. For each node in the tree, 
the number of operations to compute $U$ and
\Revector and the storage required for these vectors are both constant
times the number of beads in the node. Therefore, approximately $\bigO(n)$
operations and storage are required for each level of the tree
structure and the overall complexity and memory requirement for the
algorithm are $\bigO(n \log n)$.

\section{Hierarchical $\bigO(n)$ algorithms for $Q \vecv$ and $Q^T \vecv$} 
\label{sec:upwarddown}
In the Brownian dynamics applications, one only needs the results of
$Q\vecv$ and $Q^T\vecv$ instead of generating $Q$ and $Q^T$
explicitly. In this section, we show how to apply the hierarchical
modeling technique to further compress the information and reduce the
operations and storage for each tree node to a constant, independent
of the number of beads contained in the bead. As a result, the
overall algorithm complexity and storage both become asymptotically optimal
$\bigO(n)$. 

\subsection{Upward pass for computing $Q \vecv$}
\label{sec:upward}
We first consider $Q \vecv$ and discuss how to compress the
information in the vectors $U$ and \Revector. 
We introduce the following definitions.
\begin{definition}
For a non-leaf node $p$ containing $n>1$ beads in the tree structure,
its \Infoset \ $\bM_p$ and \Residue \  
$\vecres_p$ are respectively defined as 
\begin{equation}
\label{eq:infoset}
\bM_p = Q_{Z_p} \cdot 
\begin{bmatrix} \vecf_1 \\ \vdots \\ \vecf_n \end{bmatrix}, \quad 
\vecres_p = \Revector \cdot 
\begin{bmatrix} \vecf_1 \\ \vdots \\ \vecf_n \end{bmatrix}, 
\end{equation}
where $\vecf_j$ is the force acting on bead $j$, $Q_{Z_p}$ is the
orthogonal matrix in the decomposition of $p$'s \Zmatrix, and we
assume the origin is located at the centroid of $p$.  
For a leaf node with only one bead, we define $\bM_p = I \cdot \vecf$.
\end{definition}

Ordering $\vecf$ for parent $p$ as $[\vecf_x; \vecf_y]$, from 
Eqs.~(\ref{eq:Ufull}) and (\ref{eq:residuefull}), we have 
\begin{align*}
& \begin{bmatrix} Q_{Z_p} \\ \Revector \end{bmatrix} \cdot \vecf
= \begin{bmatrix} 
	I_{n_p} \vecf \\ U \vecf \\ \Revector \cdot 
  \vecf \end{bmatrix} \\
= & \begin{bmatrix} 
\sqrt{\frac{n_x}{n}} I_{n_x} \vecf_x + 
\sqrt{\frac{n_y}{n}} I_{n_y} \vecf_y \\
\begin{bmatrix} 
Q_{11}^T & Q_{12}^T & Q_{13}^T\\
Q_{21}^T & Q_{22}^T & Q_{23}^T\\
Q_{31}^T & Q_{32}^T & Q_{33}^T
\end{bmatrix} 
\begin{bmatrix} 
V & \veczero \\
\veczero & W \\
\sqrt{\frac{n_y}{n}} I_{n_x}  & - \sqrt{\frac{n_x}{n}} I_{n_y}
\end{bmatrix} 
\begin{bmatrix} 
\vecf_x \\ \vecf_y 
\end{bmatrix} 
\end{bmatrix} \\
= & \begin{bmatrix} 
\sqrt{\frac{n_x}{n}} \bM_x(1:3) + 
\sqrt{\frac{n_y}{n}} \bM_y(1:3) \\
\begin{bmatrix} 
Q_{11}^T & Q_{12}^T & Q_{13}^T\\
Q_{21}^T & Q_{22}^T & Q_{23}^T\\
Q_{31}^T & Q_{32}^T & Q_{33}^T
\end{bmatrix} 
\begin{bmatrix} 
\bM_x(4:6) \\ \bM_y(4:6) \\ 
\sqrt{\frac{n_y}{n}} \bM_x(1:3) - \sqrt{\frac{n_x}{n}} \bM_y(1:3) 
\end{bmatrix} 
\end{bmatrix}. 
\end{align*}
Combined with Eq.~(\ref{eq:Tformula}) that states both $T_{22}$ and
$Q_{ij}$ can be computed from children's $R_{21}$, $R_{22}$, $S_{21}$,
and $S_{22}$, all are $3\times 3$ matrices, we see that both $\bM_p$
and $\vecres_p$ can be computed using only a constant number
of operations by
\begin{align}
\bM_p &= \begin{bmatrix} 
\sqrt{\frac{n_x}{n}} \bM_x(1:3) + \sqrt{\frac{n_y}{n}} \bM_y(1:3) \\
\begin{bmatrix} Q_{11}^T & Q_{12}^T & Q_{13}^T \end{bmatrix} 
\begin{bmatrix} 
\bM_x(4:6) \\ \bM_y(4:6) \\
\sqrt{\frac{n_y}{n}} \bM_x(1:3) - 
\sqrt{\frac{n_x}{n}} \bM_y(1:3) 
\end{bmatrix} 
\end{bmatrix},  
\label{eq:implicit_MP}
\end{align}
\begin{align}
\vecres_p &= \begin{bmatrix} 
Q_{21}^T & Q_{22}^T & Q_{23}^T \\ 
Q_{31}^T & Q_{32}^T & Q_{33}^T 
\end{bmatrix} 
\begin{bmatrix} 
\bM_x(4:6) \\ \bM_y(4:6) \\ 
\sqrt{\frac{n_y}{n}} \bM_x(1:3) - 
\sqrt{\frac{n_x}{n}} \bM_y(1:3) 
\end{bmatrix}.
\label{eq:implicit_formula}
\end{align}
 
Clearly, the \Infoset of the child node contains the compressed
information for generating  parent node $p$'s \Infoset and \Residue. 
It is interesting to compare the \Infoset \ with 
the ``multipole expansion" in the fast multipole method (FMM). Both
provide effective ways to compress data contained in the node which
will be sent to the interacting nodes in the tree structure.

{\noindent \bf Bead \& Bead and Rigid Body \& Bead Cases:} 
When one or both of the child nodes are childless, the \Infoset and
\Residue of the parent can be constructed using
Eqs.~(\ref{eq:Ubeadbead}), (\ref{eq:Ubeadrigid}), and 
(\ref{eq:residuebeadrigid}). An alternative approach is to store a
size $6$ vector for both the leaf and nonleaf nodes such that a
unified formula can be used for all cases. For a leaf node, the
\Infoset is simply $\bM = [\vecf; \veczero]$, where the last $3$
numbers are due to the fact that we choose the centroid as the
origin. 

Instead of the level-wise {\tt for}-loop based execution in traditional
FMM implementation, we present in Algorithm~\ref{alg:upward} a
recursive implementation of the function {\tt Compute\_residue} for computing $Q\vecv$. 

\begin{algorithm}[htbp]
\begin{algorithmic}[1]
\Function{Compute\_residue}{$p$}
\If{node $p$ is leaf} 
\State Construct \Infoset \ $\bM_p$ directly
\Else 
\State Find child nodes $x$ and $y$ of node $p$
\State \Call{Compute\_residue}{$x$}
\State \Call{Compute\_residue}{$y$}
\State Compute $T_{22}$ and $Q_{ij}$ using Eqs.~(\ref{eq:Tformula},\ref{eq:T22full})
\State Construct $p$'s \Infoset \ $\bM_p$ and \Residue \ $\vecres$
	using Eqs.~(\ref{eq:implicit_MP},\ref{eq:implicit_formula}) 
\State Output the \Residue \ $\vecres$
\EndIf
\EndFunction
\end{algorithmic}
\caption{Recursive algorithm for computing $Q\vecv$}
\label{alg:upward} 
\end{algorithm}

{\noindent \bf Algorithm Complexity:} The complexity and storage of
the algorithm can be estimated by  
checking the operations and memory requirement for each node in the
tree structure. Notice that $T_{22}$ is 
of size $3 \times 3$, the matrix for storing $Q_{ij}$ is of size $9 \times 9$, 
\Infoset \ is a vector of size $6$, and the size of the \Residue 
 is no more than $6$. Therefore, both the number
of operations and required storage are constants for each node. The
overall algorithm complexity and  
storage are therefore proportional to the total number of nodes in the
hierarchical tree structure. For most  
practical bead distributions, as the number of nodes in the tree
structure is proportional to the number of  
beads $n$, the algorithm complexity and storage are both $\bigO(n)$.

\subsection{Downward pass for computing $Q^T \vecv$} In rigid body
dynamics, given the translational velocity $\vecV$ and angular
velocity $\vecomega$ of a rigid body $p$ consisting of $n$ beads,
located at $\{\vecr_i\}$, and
assume that the reference point for the angular velocity is located
at the centroid that is chosen as the origin, the velocity of bead $i$
in the rigid body can be computed by 
\[
\vecv_i = \vecV + \vecomega \times \vecr_i. 
\]
In matrix form, using the \Zmatrix, 
 the velocities of all the beads are given by 
\[
\vecv = Z_p^T \begin{bmatrix} \vecV \\ \vecomega \end{bmatrix}, 
\]
i.e., the velocity vector $\vecv$ is a linear combination of the column vectors
in $Z_p^T$ and the coupling coefficients are given by 
$\begin{bmatrix} \vecV \\ \vecomega \end{bmatrix}$.
To preserve orthogonality, we only consider the orthogonal basis
vectors in the decomposition of $Z_p$ and introduce the
``generalized'' velocities 
\begin{equation}
\label{eq:local}
	\left[   \begin{array}{c:c} I_{n_p}^T & U^T\end{array} \right]
  \left[ \begin{array}{c} \vecV \\ \vecomega \end{array} \right] = 
Q_{Z_p}^T 
  \left[ \begin{array}{c} \vecV \\ \vecomega \end{array} \right].
\end{equation}
Notice that only $6$ numbers in $\vecV$ and $\vecomega$ are needed to
construct the velocities of all the beads using Eq.~(\ref{eq:local}). 
We define the vector containing these numbers as the \RVset of $p$ as
follows.

\begin{definition} 
The \RVset of a rigid body $p$ consisting of $n>1$ beads 
is defined as a $6 \times 1$ vector 
$\bL = \left[ \vecV; \vecomega \right]$ 
that describes the translational and angular velocities of $p$. 
The velocities of each individual bead in the rigid body $p$ 
are given by Eq.~(\ref{eq:local}).
\end{definition}

In the hierarchical modeling technique, consider a parent node $p$
with child nodes $x$ and $y$. The velocities of all the beads of $p$
can be computed using $p$'s \RVset as $Q_{Z_p}^T \bL_p$. 
Notice that both $Q_{Z_p}$ and the \Revector computed using 
 Eqs.~(\ref{eq:Tformula}-\ref{eq:residuefull}) are combinations of 
$[Q_{Z_x}, \veczero ]$ and $[\veczero , Q_{Z_y}]$, 
which suggests that if we break the rigidity of $p$ and 
assume both $x$ and $y$ are rigid bodies but permit their relative
location to change, then both $Q_{Z_p}^T \bL_p$ 
and the linear combination of the column vectors in the transpose of 
the \Revector can be stored in 
$\left[ \begin{array}{c}    Q_{Z_x}^T\bL_x\\ Q_{Z_y}^T\bL_y \end{array}\right]$, 
 where $\bL_x$ and $\bL_y$ are the \RVset of child $x$ and $y$,
 respectively. This is summarized in the following theorem. 

\begin{theorem}
\label{theorem:downpass}
Assume a rigid body $p$ is partitioned into child $x$ and $y$ 
\[
	Q_{Z_p} = \begin{bmatrix} I_{n_p} \\ U\end{bmatrix}, \quad 
Q_{Z_x} = \begin{bmatrix} I_{n_x} \\ V \end{bmatrix}, \quad 
Q_{Z_y} = \begin{bmatrix} I_{n_y} \\ W \end{bmatrix},  
\]
then for any vector $\vecv$ containing the coupling coefficients of $p$'s \Revector, there 
exist \RVset vectors $\bL_x$ for child $x$ and $\bL_y$ for child
$y$ such that 
\begin{equation}
\label{eq:velocityeqn}
Q_{Z_p}^T \bL_p + (\Revector)^T \vecv = 
\begin{bmatrix} Q_{Z_x}^T \bL_x^T \\ Q_{Z_y}^T \bL_y^T \end{bmatrix}, 
\end{equation}
where 
\begin{align}
\bL_x &= \begin{bmatrix} 
\sqrt{\frac{n_x}{n}}I & \sqrt{\frac{n_y}{n}} Q_{13} \\ \veczero & Q_{11}
\end{bmatrix} \bL_p+
\begin{bmatrix} 
\sqrt{\frac{n_y}{n}} Q_{23} & \sqrt{\frac{n_y}{n}} Q_{33}\\
 Q_{21}  & Q_{31}
\end{bmatrix} \vecv, 
\label{eq:rvset0}
\end{align}
\begin{align}
\bL_y &= \begin{bmatrix} 
\sqrt{\frac{n_y}{n}}I & -\sqrt{\frac{n_x}{n}} Q_{13} \\
 \veczero & Q_{12} 
\end{bmatrix} \bL_p +
\begin{bmatrix} 
 -\sqrt{\frac{n_x}{ n}} Q_{23} & -\sqrt{\frac{n_x}{ n}} Q_{33} \\
  Q_{22} & Q_{32} 
\end{bmatrix}  \vecv. \label{eq:rvset1}
\end{align}
In other words, $Q_{Z_p} \bL_p + (\Revector)^T \vecv$ can be
compressed and stored in the \RVset $\bL_x$ and $\bL_y$ of child nodes
$x$ and $y$, respectively. 
\end{theorem}
Theorem~\ref{theorem:downpass} is proved by applying 
$\begin{bmatrix} Q_{Z_x} & \\ & Q_{Z_y} \end{bmatrix}$ to both sides
of Eq.~(\ref{eq:velocityeqn}) and we skip the details. 
It suggests that $Q^T\vecv$ can be computed efficiently in a downward
pass where the result of the matrix-vector product is considered as
the velocity vectors of the beads. 
Starting from the root node's \RVset, the child nodes compute their
\RVset  using Eqs.~(\ref{eq:rvset0}, \ref{eq:rvset1}) and the collection of \RVset at
leaf nodes gives the result of $Q^T \vecv$. 
When only  $\tilde{Q}^T \tilde{\vecv}$ is required, one only needs to 
replace the \RVset of the root node with a zero vector. 
For the special case when $p$ is a leaf node containing a single bead,
the velocity of the bead is $\bL_p$ as we choose the location of the
bead as the origin. 
It is also interesting to compare the \RVset \ with the ``local expansion" in 
the fast multipole methods. Both are used to store information inherited  
from parent levels. 

We present the algorithm {\tt Compute\_velocity} 
in a recursive fashion in Algorithm~\ref{alg:downward}. Here, we
assume that $Q_{ij}$ matrices for each node are already computed and
stored in the upward pass discussed in  Sec.~\ref{sec:upward}. 
\begin{algorithm}[htbp]
\begin{algorithmic}[1]
\Function{Compute\_velocity}{$p$}
\If{$p$ is a leaf node}
\State Output the \RVset \ $\bL_p$
\Else
\State Find child nodes $x$ and $y$
	\State Compute $x$ and $y$'s \RVset  using Eqs.~(\ref{eq:rvset0},\ref{eq:rvset1})
\State \Call{Compute\_velocity}{$x$}
\State \Call{Compute\_velocity}{$y$}
\EndIf
\EndFunction
\end{algorithmic}
\caption{Recursive algorithm for computing $Q^T \vecv$}
\label{alg:downward} 
\end{algorithm}

{\noindent \bf Algorithm Complexity:} Similar to the upward pass to
compute $Q \vecv$, as a constant amount of operations (and storage) is
required for each node, the complexity of the recursive algorithm is
$\bigO(n)$.

\section{Preliminary numerical results}
We present some preliminary numerical results to show the complexity
and orthogonality  properties of the new algorithms introduced in this
paper. The prototype implementation of the algorithms were done in
Matlab and the numerical tests were carried out on a personal laptop
with Intel Core i5 CPU at 2.6 GHz clock rate and 8 GB of RAM. 

\subsection{Algorithm complexity}
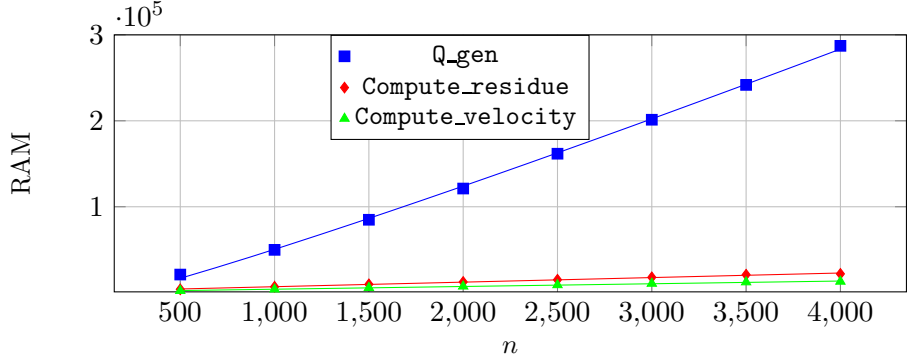
\begin{figure}[htbp]
\centering
\begin{tikzpicture} 
\begin{axis}[ 
    xlabel = {\small $n$}, 
    ylabel = {\small RAM}, 
    ymin = 1e3, 
    ymax = 3e5, 
    grid = both, 
    width = 0.8\textwidth, 
    height = 0.24\textheight, 
    legend style={at={(0.6, 1)}, font=\small}
  ] 
    \addplot [mark = square*, color = blue, only marks] coordinates {
    (500, 0.2113e5)
    (1000,0.4985e5)
    (1500, 0.8487e5)
    (2000, 1.2121e5)
    (2500, 1.6171e5)
    (3000, 2.0124e5)
    (3500, 2.4174e5)
    (4000, 2.8720e5)
  };
  \addlegendentry{{\tt Q\_{gen}}}; 
  \addplot[domain=500:4000, samples=20, forget plot, color=blue] 
  {ln(x)*x*8.87-10800}; 

  \addplot [mark=diamond*, color=red, only marks] coordinates {
    (500, 0.3948e4)
    (1000, 0.6925e4)
    (1500, 0.9776e4)
    (2000, 1.2524e4)
    (2500, 1.5115e4) 
    (3000, 1.7934e4)
    (3500, 2.0941e4)
    (4000, 2.1928e4)
  }; 
  \addlegendentry{{\tt Compute\_residue}}; 
  \addplot[domain=500:4000, samples=20, forget plot, color=red] 
  {5.31*x+1690}; 

  \addplot [mark=triangle*, color=green, only marks] coordinates {
    (500, 0.2297e4)
    (1000, 0.3968e4)
    (1500, 0.5708e4)
    (2000, 0.7405e4)
    (2500, 0.8990e4)
    (3000, 1.0635e4)
    (3500, 1.2289e4)
    (4000, 1.3038e4)
  };
  \addlegendentry{{\tt Compute\_velocity}}; 
  \addplot[domain=500:4000, samples=20, forget plot, color=green] 
  {3.17*x+908}; 
\end{axis}
\end{tikzpicture} 
\caption{Memory usage in kilobytes versus the number of beads $n$ for 
{\tt Q\_gen} (blue square), {\tt Compute\_velocity} (green triangle), 
and {\tt Compute\_residue} (red diamond).}
\label{fig:memory}
\end{figure}

We consider the memory requirement and CPU time of each algorithm for
a variety number of beads $n$ that ranges from $500$ to $4000$ with
$500$ increment. The memory requirement for each run was extracted from the
Matlab profiler's report with the memory option enabled and is
summarized in Figure~\ref{fig:memory}. Particularly, the fitted curves
are 
\begin{align*}
\mbox{\tt Q\_gen:} \quad & 8.87 n \log n - 1.08 \cdot 10^4, \\
\mbox{\tt Compute\_residue:} \quad & 5.31 n + 1.69 \cdot 10^3, \\
\mbox{\tt Compute\_velocity:} \quad & 3.17 n + 9.08 \cdot 10^2, 
\end{align*} 
which agree with the analytic results. Clearly, the explicit
generation of the matrix using {\tt Q\_gen} requires more memory than
the implicit algorithms. 

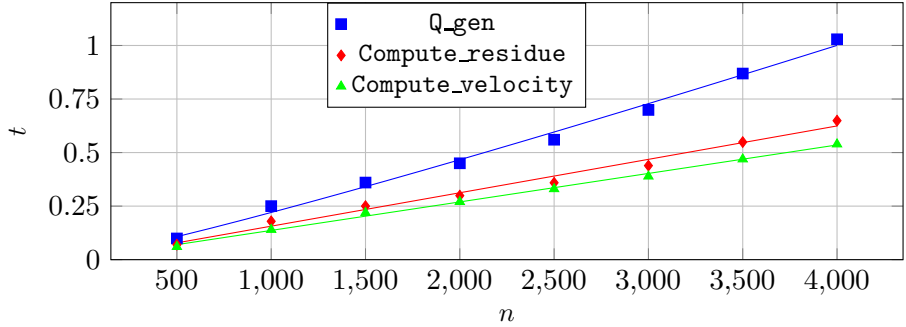
\begin{figure}[htbp]
\centering 
\begin{tikzpicture}
\begin{axis}[
    xlabel={\small $n$}, 
    ylabel={\small $t$}, 
    ytick={0,0.25,0.5,0.75,1},
    ymin = 0, 
    ymax = 1.2, 
    grid = both, 
    width = 0.8\textwidth, 
    height = 0.24\textheight, 
    legend style={at={(0.6, 1)}, font=\small}
  ]
  \addplot [mark = square*, color = blue, only marks] coordinates {
    (500, 0.099)
    (1000, 0.25)
    (1500, 0.36)
    (2000, 0.45)
    (2500, 0.56)
    (3000, 0.699)
    (3500, 0.869)
    (4000, 1.029)
  };
  \addlegendentry{{\tt Q\_{gen}}}; 
  \addplot[domain=500:4000, samples=20, forget plot, color=blue] 
  {ln(x)*x*2.97e-5+1.51e-2}; 

  \addplot [mark=diamond*, color=red, only marks] coordinates {
    (500, 0.07)
    (1000, 0.179)
    (1500, 0.25)
    (2000, 0.299)
    (2500, 0.359) 
    (3000, 0.439)
    (3500, 0.549)
    (4000, 0.649)
  }; 
  \addlegendentry{{\tt Compute\_residue}}; 
  \addplot[domain=500:4000, samples =20, forget plot, color=red] 
  {1.56e-4*x + 3.57e-4}; 

  \addplot [mark=triangle*, color=green, only marks] coordinates {
    (500, 0.06)
    (1000, 0.14)
    (1500, 0.219)
    (2000, 0.27)
    (2500, 0.33)
    (3000, 0.389)
    (3500, 0.469)
    (4000, 0.539)
  };
  \addlegendentry{{\tt Compute\_velocity}}; 
  \addplot[domain=500:4000, samples=20, forget plot, color=green] 
  {1.33e-4*x + 3.57e-3}; 

 \end{axis} 
\end{tikzpicture} 
\caption{CPU time in seconds ($t$) versus the number of beads $n$ for 
{\tt Q\_gen} (blue square), {\tt Compute\_velocity} (green triangle), 
and {\tt Compute\_residue} (red diamond).} 
\label{fig:cputime}
\end{figure}

In Figure \ref{fig:cputime}, we present the CPU time in seconds for
the $3$ algorithms. In the experiments, each algorithm was executed
ten times for each value of $n$ and we present $min_{ 1 \leq i \leq 10} (\mbox{CPU time in $i^{th}$ run})$.
The fitted curves are 
\begin{align*}
\mbox{\tt Q\_gen:} \quad & 2.97\cdot 10^{-5} n \log n + 1.51\cdot 10^{-2}, \\
\mbox{\tt Compute\_residue:} \quad & 1.56\cdot 10^{-4} n + 3.57 \cdot 10^{-4}, \\
\mbox{\tt Compute\_velocity:} \quad & 1.33\cdot 10^{-4} n + 2.57 \cdot 10^{-3}, 
\end{align*} 
and match the analytic results. The implicit methods are more
efficient than the explicit function {\tt Q\_gen}. Furthermore, the
downward pass {\tt  Compute\_velocity} is more efficient than the
upward pass {\tt Compute\_residue} because the upward pass needs to
compute the matrices $T_{22}$ and $Q_{ij}$ that are also used in the
downward pass. 

\subsection{Orthogonality properties}
By introducing the orthogonal numerical linear algebra technique, the
hierarchical models and recursive algorithms better preserve the
orthogonality properties of the matrix $Q$ and
$Q^T$. Table~\ref{tab:result} shows various orthogonality properties
of the algorithms. All the errors are measured in infinity norm. In
columns two and three, we consider $QQ^T- I$ and $Q^TQ-I$ of matrix
$Q$ explicitly generated using {\tt Q\_gen}. In column four, we
apply the explicitly generated $Q$ to a given vector $\vecv$ and
compare the output to that computed from {\tt Compute\_residue}. In
column five, we apply the transpose of the explicitly generated matrix
$Q$ to a given vector and compare the output to that computed from
{\tt Compute\_velocity}. In column six, we first apply {\tt
  Compute\_residue} to compute $Q\vecv$ and then apply function 
{\tt Compute\_velocity} and compare it against the input vector. In
column seven, we first apply {\tt Compute\_velocity} to compute $Q^T
\vecv$ and then apply function {\tt Compute\_residue} and compare it
against the input vector. All the errors are close to machine
precision, providing strong evidence of the orthogonality preserving
properties of the three algorithms. 
 
\begin{table}
\centering
\begin{tabular}{|c |c|c|c|c|c|c|} 
\hline 
$n$ & $QQ^T - I$ & $Q^TQ-I$ & $Q\vecv$ & $Q^T\vecv$ & 
$Q^T(Q\vecv) - \vecv$ & $Q(Q^T\vecv) - \vecv$ \\
\hline 
500 & 3.1e-15 & 6.6e-16 & 2.8e-15 & 2.3e-15 & 1.1e-15 & 4.8e-15 \\ \hline
1000 & 1.3e-15 & 8.8e-16 & 8.4e-15 & 2.4e-15 & 2.2e-15 & 9.7e-15\\ \hline
1500 & 2.2e-14 & 1.1e-15 & 1.4e-14 & 1.8e-15 & 3.0e-15 & 7.3e-15\\ \hline
2000 & 1.9e-15 & 1.9e-15 & 4.0e-14 & 3.1e-15 & 1.8e-15 & 1.3e-14\\ \hline
2500 & 6.2e-15 & 8.8e-16 & 2.0e-14 & 1.8e-15 &  2.0e-15 & 1.7e-14\\ \hline
3000 & 7.1e-15 & 1.1e-15 & 1.9e-14 & 4.5e-15 & 3.2e-15 & 1.9e-14\\ \hline
3500 & 7.4e-14 & 1.1e-15 & 3.0e-14 & 5.6e-15 & 1.7e-15 & 2.8e-14\\ \hline
4000 & 2.8e-15 & 1.1e-15 & 3.4e-14 & 7.2e-15 & 1.6e-15 & 2.1e-14\\ \hline
8000 & 1.2e-13 & 8.8e-16 & 9.4e-14 & 5.3e-15 & 3.6e-15 & 2.3e-14\\ \hline
10000 & 7.1e-15 & 8.8e-16 & 7.5e-14 & 1.1e-14 & 2.9e-15 & 3.5e-14\\ \hline
\end{tabular}  
\caption{Orthogonality preserving quality of Algorithms 1-3.}
\label{tab:result} 
\end{table}

\section{Summary and future work}
In this paper, we apply the hierarchical modeling technique and
present $3$ recursive algorithms for generating special hierarchical
structured orthogonal matrices with applications in Brownian dynamics
simulations of biomolecular systems. By combining the orthogonal
linear algebra techniques with the hierarchical models, our
preliminary numerical experiments show that the implicit algorithms
for computing $Q \vecv$ and $Q^T \vecv$ are both asymptotically
optimal in complexity and have good orthogonality properties.

We are currently implementing the parallel versions of
Algorithms~\ref{alg:upward} and \ref{alg:downward} and developing
toolboxes that will be integrated with our Brownian dynamics
simulations package. As the algorithmic structure of these models is
very close to the fast multipole method, we plan to adapt and extend
parallelization techniques in existing packages such as 
{\sc rec}FMM \cite{zhang2016recfmm} and DASHMM \cite{zhang2016dashmm}.
 The developed software will be released to the research community
 under open-source license agreement.

Finally, we want to mention that other reformulations of Eq.~(\ref{eq:rigidbodyBDHI}) 
are also possible. For instance,  the Schur complement was used for a similar 
problem\cite{usabiaga2016hydrodynamics}.  Comparisons of different
reformulations are being performed by the authors. Another closely
related research topic is the design of  effective preconditioners for 
Eq.~(\ref{eq:newBDHI}). Some of these topics are briefly discussed in \cite{fang2016recursive}. 
Detailed results along these directions will be discussed in future papers. 

\bibliographystyle{plain} 
\bibliography{rigidbody,huangbib}

\end{document}